\documentclass[final,reqno,12pt]{amsart}
\usepackage{amsmath}
\usepackage{amssymb}
\usepackage{enumerate}
\usepackage{caption}
\usepackage{url}
\usepackage[colorlinks=true, citecolor=red, urlcolor=blue]{hyperref}\newcommand{\HLink}[2]{\href{#1}{\underline{#2}}}
\usepackage[pdftex,final]{graphicx}  

\newcommand{\SpecialText}[1]{{\text{#1}}}
\newcommand{\SpecialTextFN}[1]{{\text{#1}}}

\def\medskip{\par\vspace{10pt}}
\def\smallskip{\par\vspace{6pt}}

\newcommand{\FigureCaption}[1]{%
\captionsetup{labelfont={small,bf},textfont={small,rm},labelsep=period,singlelinecheck=false}
\caption{#1}}

\newcommand{\webpage}[1]{Available at \[\hbox{\url{#1}}\]}
\newcommand{\webpageinline}[1]{Available at \url{#1}}

\newcommand{\Riemann}[1]{\mathcal{#1}}
\newcommand{\RiemannMap}[3]{\mathrm{#1}\xrightarrow{\Riemann{#2}}\mathrm{#3}}
\newcommand{\RiemannMapB}[2]{\xrightarrow{\Riemann{#1}}\mathrm{#2}}
\newcommand{\Tonn}{\mathcal{T}}
\newcommand{\TonnMap}[2]{\mathrm{#1}\xrightarrow{\Tonn}\mathrm{#2}}
\newcommand{\TonnMapB}[1]{\xrightarrow{\Tonn}\mathrm{#1}}
\newcommand{\TonnMapTT}[2]{\mathrm{#1}\xrightarrow{\Tonn\Tonn}\mathrm{#2}}
\newcommand{\TonnMapTTB}[1]{\xrightarrow{\Tonn\Tonn}\mathrm{#1}}
\newcommand{\Tonnetz}{\emph{Tonnetz}}

\newcommand{\Emph}[1]{\textbf{\emph{#1}}}

\newcommand{\Item}[1]{\item[]}

\newcommand{\IML}{\emph{In My Life}}
\newcommand{\Chordino}{\textsc{Chordino}}
\newcommand{\Audacity}{\textsc{Audacity}}
\title{A Geometric Analysis Of The Harmonic Structure of \emph{In My Life}}

\author{James S.~Walker}
\address{James S.~Walker\\
         Mathematics\\
         Univ.~of Wisconsin-Eau Claire\\
         walkerjs@uwec.edu}

\author{Gary W.~Don}
\address{Gary W.~Don\\
         Music \& Theatre Arts\\
         Univ.~of Wisconsin-Eau Claire\\
         dongw@uwec.edu}

\begin{document}
 \ \vspace{-80pt}
\begin{abstract}
After our book \cite{book:WalkerAndDon} was published, we found a striking example of the importance of
the \Tonnetz\ for analyzing the harmonic structure of The Beatles' song,
\IML. Our \Tonnetz\ analysis will illustrate the highly structured geometric
logic underlying the numerous chord progressions in the song. Spectrograms
provide a way for us to visualize chordal harmonics and their connection with
voice leading. We shall also describe the interesting harmonic rhythms of the
song's chord progressions. A lot of this harmonic rhythm lends itself well
to a geometric description.
\end{abstract}

\maketitle

\section{Introduction\label{sec:Intro}}%
The Beatles, especially in some classic songs of John Lennon and Paul McCartney, used much more elaborate chord
progressions than one finds in typical popular songs. A fine example of this occurs with their
song, \IML. For whatever reason, the chord progressions and instrumental melodies found in various sheet musics
for Lennon and McCartney songs, such as \IML, are aural transcriptions made by
others.\footnote{We will not speculate whether this is due to the common assertion that Lennon and McCartney could
not read music, or due to them simply leaving musical transcriptions of their songs to others.}
This has led to some differences in the various sheet musics, and chord listings, that are available. Our approach to this difficulty is to do an initial analysis on a basic template for the song found in the lead sheet version
in \cite{url:LeadSheet}.
After analyzing its harmonic structure, consisting of a \Tonnetz\ analysis of its chord progressions and an examination
of its harmonic rhythm, we then look at a computer analysis of the recorded performance of \IML\ from
The Beatles' album \emph{Rubber Soul}. The harmonic structure of this recording involves just a few enhancements
of the lead sheet version.
\section{Harmonic Structure of the Lead Sheet Version of \IML\label{sec:HarmonicStructureLeadSheet}}%
The lead sheet lays out the song form for \IML, in the key of $\mathrm{A}$. We show this song form here:
\begin{equation}\label{eq:SongForm}
\SpecialText{Verse}\quad\ \SpecialText{Bridge}\quad\ \SpecialText{Verse}\quad\ \SpecialText{Bridge}\quad\ \SpecialText{Interlude}
\quad\ \SpecialText{Bridge}\quad\ \SpecialText{Coda}
\end{equation}
In this section, we shall analyze the harmonic structure of each distinct part---$\SpecialText{Verse}$, $\SpecialText{Bridge}$,
$\SpecialText{Interlude}$, and $\SpecialText{Coda}$---of this song form.
\subsection{Harmonic Structure of $\SpecialText{Verse}$\label{subsec:HarmonicStructureVerse}}%
The lead sheet in \cite{url:LeadSheet} gives the chord progressions for $\SpecialText{Verse}$. These progressions, in
terms of \Tonnetz\ transformations, are
\begin{equation}\label{eq:VerseProgressions}
\TonnMap{A}{E^7}\TonnMapB{A}\TonnMapB{E^7}\TonnMapB{A}
\TonnMapB{f^\sharp}\TonnMapB{A^{\! 7}}\TonnMapB{D}\TonnMapB{d}\TonnMapB{A}
\TonnMapB{f^\sharp}\TonnMapB{A^{\! 7}}\TonnMapB{D}\TonnMapB{d}\TonnMapB{A}
\end{equation}
In the second instance of $\SpecialTextFN{Verse}$, the progression $\TonnMap{A}{E^7}$ only occurs once.
It is worth noting that all of these progressions are single \Tonnetz\ transformations.
They follow an interesting path through chords on the \Tonnetz.

\begin{figure}[!htb]
\setlength{\unitlength}{1in}
\begin{picture}(6.5,3.4)\thicklines
\put(0.6,3.3){\SpecialText{Verse}}
\put(0,0){\resizebox{1.5in}{!}{\includegraphics{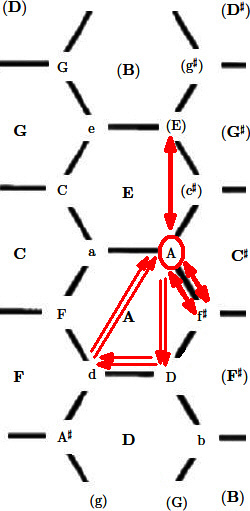}}}
\put(1.6,0){\color{green}\line(0,1){3.4}}
\put(2.3, 3.3){\SpecialText{Bridge}}
\put(1.7,0){\resizebox{1.5in}{!}{\includegraphics{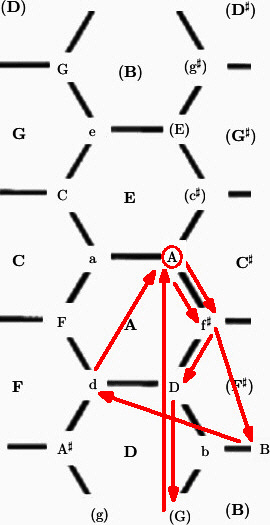}}}
\put(3.25,0){\color{green}\line(0,1){3.4}}
\put(3.8,3.3){\SpecialText{Interlude}}
\put(3.3,0){\resizebox{1.5in}{!}{\includegraphics{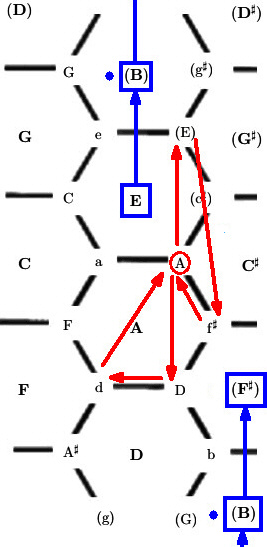}}}
\put(4.9,0){\color{green}\line(0,1){3.4}}
\put(5.6,3.3){\SpecialText{Coda}}
\put(5,0){\resizebox{1.5in}{!}{\includegraphics{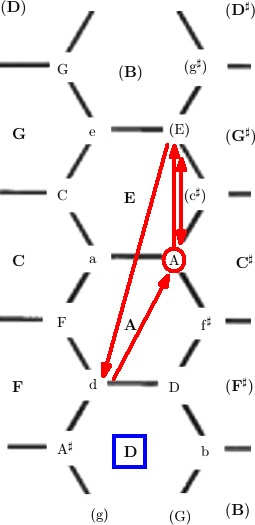}}}
\end{picture}
\FigureCaption{Red arrows show chord progressions on the \Tonnetz\ for each of the parts in the
song. In $\SpecialTextFN{Verse}$, $\SpecialTextFN{Interlude}$, and $\SpecialTextFN{Coda}$, the circled chord $\mathrm{A}$ is the starting and ending chord. In $\SpecialTextFN{Bridge}$, the circled chord $\mathrm{A}$ is the ending chord of $\SpecialTextFN{Verse}$ or $\SpecialTextFN{Interlude}$, hence $\TonnMap{A}{f^\sharp}$ is the transition to the beginning of $\SpecialTextFN{Bridge}$. The chord $\mathrm{A}$ is also the ending chord of $\SpecialTextFN{Bridge}$. The blue rectangles, and connecting arrows, in $\SpecialTextFN{Interlude}$ illustrate the movement of pitches that undergirds the double \Tonnetz\ transformation $\TonnMapTT{E}{f^\sharp}$. The blue dot next to the pitch class $\mathrm{B}$ indicates its double position in the \Tonnetz. The blue rectangle in $\SpecialTextFN{Coda}$ marks the note $\mathrm{D}$ found in both the $\mathrm{E}^7$ chord and the $\mathrm{d}$ chord.\label{fig:TonnetzDiagramsForInMyLife}}
\end{figure}

In Figure~\ref{fig:TonnetzDiagramsForInMyLife}, we have plotted the motion of these progressions through the
chords on the \Tonnetz. In these plots, any seventh chords are treated as embellishments of underlying triadic chords.
For $\SpecialText{Verse}$, the motions of the chord progressions are shown on the left side of Figure~\ref{fig:TonnetzDiagramsForInMyLife}.

\begin{figure}[!htb]\centering
\setlength{\unitlength}{1in}
\begin{picture}(3,3)
\put(0,0){\resizebox{!}{3in}{\includegraphics{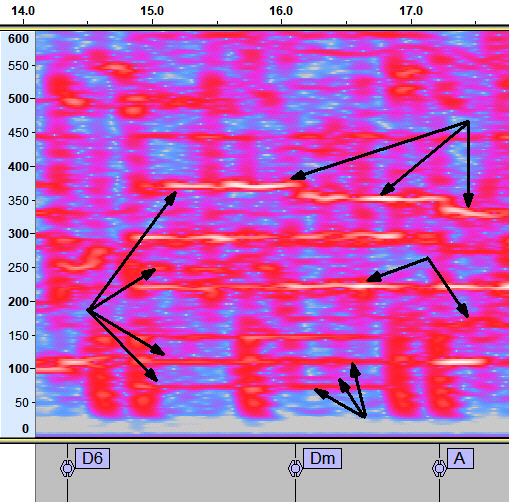}}}
\end{picture}
\FigureCaption{The quadruple arrow points to harmonics for pitches of type $\mathrm{D}, \mathrm{A}, \mathrm{B}, \mathrm{F}^\sharp$ in ascending order. This is why computer analysis classifies the chord as $\mathrm{D}^6$ (more on
this computer analysis in Section~\ref{sec:HarmonicStructureRecording}). In the \Tonnetz\ we have plotted this as a $\mathrm{D}$ chord, as indicated in the lead sheet. The triple arrow at the bottom points to fundamentals for a $\mathrm{d}$-chord
($\mathrm{D}_2, \mathrm{F}_2, \mathrm{A}_2$). The double arrow points to fundamentals for the chromatic descent from an $\mathrm{A}_3$ note to a natural $\mathrm{F}_3$ note, and the fundamental of $\mathrm{F}_3$ aligns with a previous second harmonic for the $\mathrm{F}_2$ note in the $\mathrm{d}$-chord. The triple arrow at the top points to harmonics for the voice leading
$\mathrm{F^\sharp}\to \mathbf{F^\natural}\to\mathbf{E}$. \label{fig:VoiceLeadingsAndHarmonicsInVerseClip}}
\end{figure}

The ending progression $\TonnMap{A^{\! 7}}{D}\TonnMapB{d}\TonnMapB{A}$, corresponding to the triangular path around
the pitch class hexagon $\mathbf{A}$ is particularly interesting. Here $\mathrm{A^{\! 7}}$ is functioning
as a \emph{secondary dominant} for the chord $\mathrm{D}$, i.e., a dominant for the non-tonic chord $\mathrm{D}$. In roman
numerals, this is denoted as $\mathbf{V}^7 / \mathbf{IV}$. Hence, $\TonnMap{\mathbf{V}^7 / \mathbf{IV}}{\mathbf{IV}}$ is the roman numeral version of $\TonnMap{A^{\! 7}}{D}$ in the key of $\mathrm{A}$. The transformation $\TonnMap{D}{d}$ is
a neo-Riemannian transformation $\RiemannMap{D}{P}{d}$ which introduces a \emph{modal mixture}. It is well-known
that The Beatles were experimenting with songs in various modes and modal mixtures on the \emph{Rubber Soul} album,
and the \emph{Revolver} album that followed it. For example, the Mixolydian mode is used in the song \emph{Norwegian Wood}
on \emph{Rubber Soul}. There is also a mode mixture used in \emph{Eleanor Rigby} on \emph{Revolver}, which we
discussed in some detail in \cite[Example~7.4.6]{book:WalkerAndDon}. For \IML, using the minor chord $\mathrm{d}$---which
belongs to the key of $\mathrm{a}$-minor parallel to the key of $\mathrm{A}$-major---provides a mode mixture.
The particular progression of $\RiemannMap{\mathbf{IV}}{P}{\mathbf{iv}}$ is often used in music to provide a nice contrast
in ``color'' of these parallel chord types. Within the song \IML, $\RiemannMap{D}{P}{d}$ harmonizes
well with the chromatic descent in the lyrics from an $\mathrm{A}$ note to a natural $\mathrm{F}$ note. See
Figure~\ref{fig:VoiceLeadingsAndHarmonicsInVerseClip}. Moreover, when the ending chord $\mathrm{A}$ is included to form
the cadence $\RiemannMap{D}{P}{d}\RiemannMapB{N}{A}$, there is a smooth voice leading between the pitches $\mathbf{F^\sharp}\to \mathbf{F^\natural}\to\mathbf{E}$ in the respective chords. See again, Figure~\ref{fig:VoiceLeadingsAndHarmonicsInVerseClip}.
In roman numerals, this cadence is $\mathbf{IV}\to\mathbf{iv}\to\mathbf{I}$, i.e., a plagal cadence with modal mixture. We have focused here on the ending cadence, but the other progressions and their motions on
the \Tonnetz\ are also worth considering, especially in connection with the harmonic rhythm in $\SpecialText{Verse}$.

The harmonic rhythm in $\SpecialText{Verse}$ closely relates to the structure of the path followed by the
chord changes on the \Tonnetz. This path breaks down into three substructures:
\begin{equation}\label{eq:SubstructuresForVerseProgressions}
\hbox{%
\begin{tabular}{ccccc}
$\TonnMap{A}{E^7}$ &\qquad\qquad
& $\RiemannMap{A}{R}{f^\sharp}\RiemannMapB{R}{A^{\! 7}}$
&\qquad\qquad & $\RiemannMap{D}{P}{d}\RiemannMapB{N}{A}$ \\[5pt]
$\TonnMap{\mathbf{I}}{\mathbf{V^7}}$ &\qquad\qquad
&$\RiemannMap{\mathbf{I}}{R}{\mathbf{vi}}\RiemannMapB{R}{\mathbf{V^7/IV}}$
&\qquad\qquad & $\RiemannMap{\mathbf{IV}}{P}{\mathbf{iv}}\RiemannMapB{N}{\mathbf{I}}$
\end{tabular}}
\end{equation}
Each of these substructures is repeated twice in moving through the \Tonnetz.
The roman numerals are for the chords in the key of $\mathrm{A}$-major, including the
notation $\mathbf{V^7/IV}$ that we previously discussed.

\begin{figure}[!htb]\centering
\setlength{\unitlength}{1.1in}
\begin{picture}(6.5,1.5)
\put(0.5,0){
\begin{picture}(6.5,1.5)
\put(0.15,0){
\begin{picture}(1.5,1.5)
\put(0.56,1.42){$\mathrm{A}$}\put(0.53,-0.03){$\mathrm{E^7}$}
\put(0,0.15){\resizebox{1.32in}{!}{\includegraphics{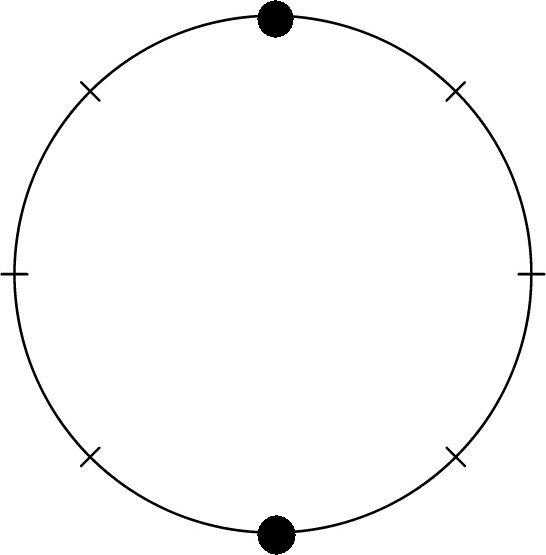}}}
\put(0.4,0.71){{\tiny $\TonnMap{\mathbf{I}}{\mathbf{V}^7}$}}
\end{picture}}
\put(1.4,1.0){\vector(1,0){0.3}}
\put(1.6,0){
\begin{picture}(1.5,1.5)
\put(0.76,1.42){$\mathrm{A}$}\put(0.76,-0.03){$\mathrm{f^\sharp}$}\put(0.01,0.71){$\mathrm{A^{\! 7}}$}
\put(0.2,0.15){\resizebox{1.32in}{!}{\includegraphics{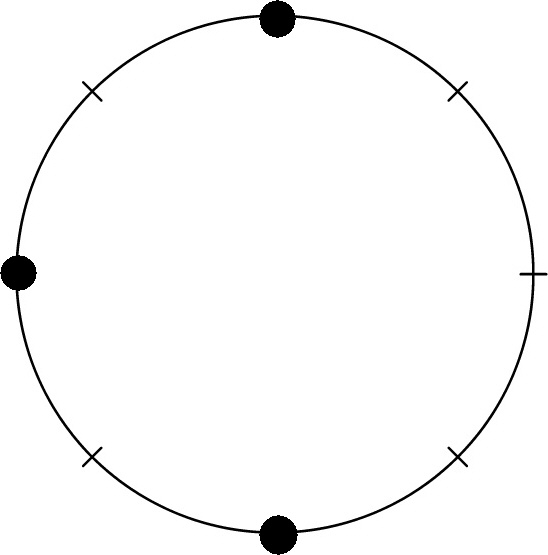}}}
\put(0.37,0.71){{\tiny $\RiemannMap{\mathbf{I}}{R}{\mathbf{vi}}\RiemannMapB{R}{\mathbf{V^7/IV}}$}}
\end{picture}}
\put(3.1,1.0){\vector(1,0){0.3}}\put(3.3,1.0){\vector(-1,0){0.2}}
\put(3.2,0){
\begin{picture}(1.5,1.5)
\put(0.75,1.42){$\mathrm{D}$}\put(0.75,-0.03){$\mathrm{A}$}\put(1.45,0.71){$\mathrm{d}$}
\put(0.2,0.15){\resizebox{1.32in}{!}{\includegraphics{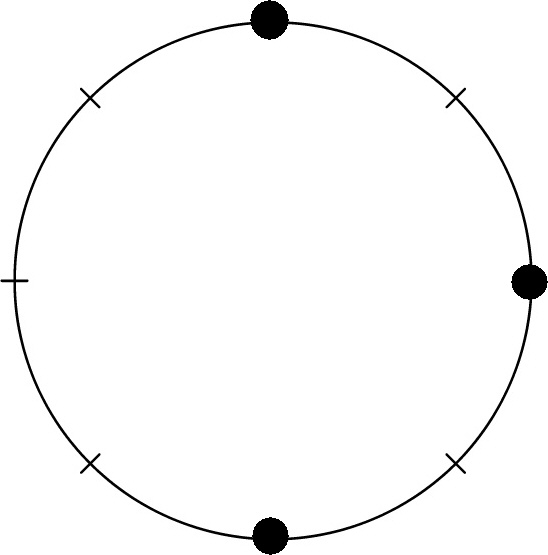}}}
\put(0.445,0.71){{\tiny $\RiemannMap{\mathbf{IV}}{P}{\mathbf{iv}}\RiemannMapB{N}{\mathbf{I}}$}}
\end{picture}}
\end{picture}}
\end{picture}
\FigureCaption{Harmonic rhythm of $\SpecialTextFN{Verse}$. The eight hour clocks mark $1$ hour per quarter note beat, spanning two measures each. In the first instance of $\SpecialTextFN{Verse}$, the clock on the left is cycled through twice. While in the second instance, it is cycled through once. In both instances, the second and third clocks are each cycled through twice in alternation. The arrows mark the change from one clock to the next---one harmonic rhythm to the next---as the music proceeds through $\SpecialTextFN{Verse}$, two measures at a time. Within each clock, there is a labeling of its form in terms of roman numerals and types of transformations (either general \Tonnetz\ transformation $\Tonn$, or neo-Riemannian transformation $\Riemann{R}, \Riemann{P}, \text{or}\ \Riemann{N}$).\label{fig:HarmonicRhythmOfVerse}}
\end{figure}

The substructures in \eqref{eq:SubstructuresForVerseProgressions} correspond to three distinct harmonic rhythms.
In Figure~\ref{fig:HarmonicRhythmOfVerse}, we show these three distinct harmonic rhythms by plotting them on
rhythm clocks of eight hours each, one hour per beat (quarter note), spanning two measures each. The first substructure
in \eqref{eq:SubstructuresForVerseProgressions} has a whole note rhythm (one chord per measure), while the
second and third substructures show alternating faster (half-note) rhythm with slower (whole note) rhythm. Furthermore,
it is interesting that the second and third clocks in Figure~\ref{fig:HarmonicRhythmOfVerse} are reflections of each
other. The reflections being through a mirror passing through hours $0$ and $4$, which are the hours for the chord
onsets in the first clock.

\subsection{Harmonic Structure of $\SpecialText{Bridge}$\label{subsec:HarmonicStructureBridge}}%
The lead sheet in \cite{url:LeadSheet} gives the following chord progressions for
$\SpecialText{Bridge}$:
\begin{equation}\label{eq:BridgeProgressions}
(\TonnMap{A}{}\!\!)\;\TonnMap{f^\sharp}{D}\TonnMapB{G}\TonnMapTTB{A}
\TonnMapB{f^\sharp}\TonnMapB{B^7}\TonnMapTTB{d}\TonnMapB{A}
\end{equation}
The progression $(\TonnMap{A}{}\!\!)$ is not part of $\SpecialTextFN{Bridge}$.
It indicates the progression from $\SpecialTextFN{Verse}$ to
$\SpecialTextFN{Bridge}$, or $\SpecialTextFN{Interlude}$ to $\SpecialTextFN{Bridge}$.

These transformations are graphed on the \Tonnetz\ in Figure~\ref{fig:TonnetzDiagramsForInMyLife}.
Most of them are single \Tonnetz\ transformations. These single \Tonnetz\ transformations
either pass across a single hexagon from the set $\{\mathbf{D}, \mathbf{F}^\sharp, \mathbf{A}\}$,
corresponding to the chord $\mathrm{D}$, or along an edge between two of these pitch
classes corresponding to $\mathrm{D}$. In fact, the $\mathrm{D}$-chord lies in the center
of the network of transformations shown for $\SpecialText{Bridge}$ in the \Tonnetz\ diagram.
All these geometric facts provide clear evidence for the \emph{tonicization}\footnote{Tonicization refers to a strong
emphasis, or central role, for a scale degree or chord other than the tonic for the key.} of the subdominant
chord $\mathrm{D}$ in $\SpecialText{Bridge}$.
\begin{figure}[!htb]\centering
\setlength{\unitlength}{1in}
\begin{picture}(2.2,3)
\put(0,0){\resizebox{!}{3in}{\includegraphics{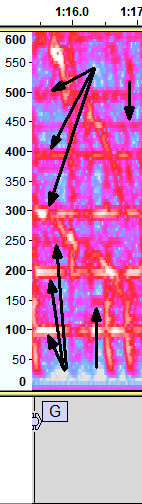}}}
\put(1.5,0){\resizebox{!}{3in}{\includegraphics{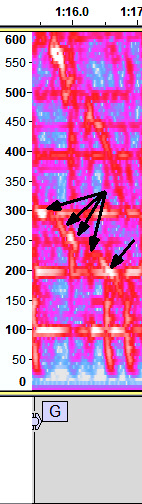}}}
\end{picture}
\FigureCaption{Left: The triple arrow at the bottom indicates harmonics
with frequencies corresponding to pitches $\mathrm{G}_2, \mathrm{G}_3, \mathrm{B}_3$ in ascending order. While the triple
arrow at the top indicates harmonics with frequencies corresponding to pitches $\mathrm{D}_4, \mathrm{G}_4, \mathrm{B}_5$ in ascending order. Consequently, computer analysis classifies the chord being played as a $\mathrm{G}$-chord. The downward
pointing arrow points to a harmonic for $\mathrm{A}_4$ which is an overtone for the harmonic for $\mathrm{D}_4$ pointed to
by the upward pointing arrow. Right: The quadruple arrow points to fundamentals for the descending sequence of notes
$\mathrm{D}, \mathrm{C^\sharp}, \mathrm{B}, \mathrm{A}$. The first fundamental, for the $\mathrm{D}$ note, matches
the $\mathrm{D}_4$ harmonic in the $\mathrm{G}$-chord. The single arrow points to an implicit $\mathrm{G}$ note occurring
when Lennon's pitch descent from $\mathrm{A}$ crosses the $\mathrm{G}_3$ harmonic. A video illustrating this implicit
note is available at \cite{url:LennonLyricsOnYouTube}.\label{fig:VoiceLeadingAndHarmonicsInBridgeClip}}
\end{figure}

The two exceptions to single \Tonnetz\ transformations in \eqref{eq:BridgeProgressions} are the
double \Tonnetz\ transformations $\TonnMapTT{\mathrm{G}}{A}$ and $\TonnMapTT{\mathrm{B^7}}{d}$.
In the chord $\mathrm{G}$, the leading tone $\mathrm{G}^\sharp$ has been dropped by a half step
to $\mathrm{G}$-natural. This lowering of the leading tone by a half step in order to play a
$\mathbf{VII}$-chord was commonly employed in popular music at that time. In which case,
the roman numeral form of $\TonnMapTT{\mathrm{G}}{A}$ represents a return to the tonic with
$\mathbf{VII}\to\mathbf{I}$. The musical practice here is similar to the Renaissance practice of
lowering the seventh tone in a major scale in order to avoid tritones.
On the left of Figure~\ref{fig:VoiceLeadingAndHarmonicsInBridgeClip}, we show the harmonics for
the $\mathrm{G}$-chord. It is interesting that the notes sung by Lennon during the playing of this
chord are a descending scale sequence $\mathrm{D}$, $\mathrm{C^\sharp}$, $\mathrm{B}$, $\mathrm{A}$.
When the seventh tone for the $\mathrm{A}$-major scale is lowered by a half-step, the
resulting scale is $\mathrm{A}$-Mixolydian:
\[
\mathrm{A},\ \ \mathrm{B},\ \ \mathrm{C}^\sharp,\ \ \mathrm{D},\ \ \mathrm{E},\ \
\mathrm{F}^\sharp,\ \ \mathrm{G},\ \ \mathrm{A}.
\]
As shown on the right of Figure~\ref{fig:VoiceLeadingAndHarmonicsInBridgeClip}, Lennon's lyrical
pitch descent from $\mathrm{A}$ crosses the second harmonic of the fundamental
for the $\mathrm{G}$-chord. If we take this as an implicit sounding of a $\mathrm{G}$-note, then
the descent is along the $\mathrm{A}$-Mixolydian scale. It is interesting to note in this connection that
the note $\mathrm{G}^\sharp$ is never used in $\SpecialText{Bridge}$, neither in the lyric's notes nor in the
chords. Thus, perhaps we have at least an implicit invocation of the $\mathrm{A}$-Mixolydian mode in
$\SpecialText{Bridge}$.
As for the progression, $\TonnMapTT{\mathrm{B^7}}{d}$, when we analyze the recording of \IML, we shall
see that $\TonnMapTT{B^7}{d}$ does \Emph{not} occur there. Instead, it is replaced by $\TonnMap{B}{D}\TonnMapB{d}$.
Thus, in the recording, even though the chords $\mathrm{B}$ and $\mathrm{d}$ are not
in the key of $\mathrm{A}$, they still play a consonant role in the music due to the
\Tonnetz\ transformation sequence $\TonnMap{f^\sharp}{B}\TonnMapB{D}\TonnMapB{d}$.

There is also a change in harmonic rhythm in \SpecialText{Bridge}. The harmonic rhythm 
in \cite{url:LeadSheet} is a simple one chord per measure rhythm. This contrasts with
the slightly more complex harmonic rhythm of \SpecialText{Verse}. So we have an interesting change
of complexity of the harmonies when \SpecialText{Verse} changes to \SpecialText{Bridge}: the chord
progressions become more elaborate, while the harmonic rhythm becomes less elaborate.
That change of complexity then reverses when \SpecialText{Bridge} changes back to \SpecialText{Verse}.
\subsection{Harmonic Structure of $\SpecialText{Interlude}$\label{subsec:HarmonicStructureInterlude}}%
The lead sheet in \cite{url:LeadSheet} gives the following chord progressions for
$\SpecialText{Interlude}$:
\begin{equation}\label{eq:InterludeProgressions}
\TonnMap{A}{E}\TonnMapTTB{f^\sharp}\TonnMapB{A^{\! 7}}\TonnMapB{D}\TonnMapB{d}\TonnMapB{A}
\end{equation}
All of the transformations in \eqref{eq:InterludeProgressions} are single \Tonnetz\ transformations, except for
$\TonnMapTT{\mathrm{E}}{f^\sharp}$. An interesting feature of this one exception is that, while
the chord $\mathrm{E}$ is sounding, the notes being played are a descending series of
perfect fourths: $\mathrm{E}, \mathrm{B}, \mathrm{F^\sharp}$. We show the pitch classes for
those notes in Figure~\ref{fig:TonnetzDiagramsForInMyLife}, and it is clear that they lead nicely
into the chord $\mathrm{f^\sharp}$. It is clear from the \Tonnetz\ diagram in
Figure~\ref{fig:TonnetzDiagramsForInMyLife} that the tonic chord for $\SpecialText{Interlude}$ is
$\mathrm{A}$. This makes sense as $\SpecialText{Interlude}$ is a substitute for $\SpecialText{Verse}$
in the song form shown in \eqref{eq:SongForm}.

The harmonic rhythm in $\SpecialText{Interlude}$ is mostly half-note.
The only exception is that for each of the two passes through $\SpecialText{Interlude}$ in the lead
sheet score, the ending chord $\mathrm{A}$ is held for a full measure.
\subsection{Harmonic Structure of $\SpecialText{Coda}$\label{subsec:HarmonicStructureCoda}}%
There is a brief coda at the end of the song. The lead sheet in \cite{url:LeadSheet} gives the following
chord progressions for $\SpecialText{Coda}$:
\begin{equation}\label{eq:CodaProgressions}
\TonnMap{A}{E^7}\TonnMapTTB{d}\TonnMapB{A}\TonnMapB{E^7}\TonnMapB{A}
\end{equation}
All of these are single \Tonnetz\ transformations, with one exception: $\TonnMapTT{E^7}{d}$. This is a surprising
progression, as the dominant seventh chord $\mathrm{E^7}$ would typically anticipate the tonic chord $\mathrm{A}$.
However, the progression to the subdominant minor chord $\mathrm{d}$, used frequently throughout the song---along
with Lennon singing a chromatic descending sequence of pitches $\mathrm{C^\natural}, \mathrm{B}, \mathrm{A}$---provides
a nice unexpected twist to the music in this final part.

The harmonic rhythm for the $\SpecialText{Coda}$ is whole note throughout. It thus provides
a slowing down from the half-note harmonic rhythm at the end of the last $\SpecialText{Bridge}$.
\subsection{Summary\label{subsec:SummaryLeadSheet}}%
The lead sheet provides a foundation for understanding the fascinating harmonic structure
of \IML. We have seen that the \Tonnetz\ provides a powerful tool for understanding
the geometric logic underlying the chord progressions in the song. A geometric analysis also
sheds light on some of its harmonic rhythm.
\section{Harmonic Structure of the Recording of \IML\label{sec:HarmonicStructureRecording}}%
The recording of \IML\ for the most part follows the structure given in the lead sheet in \cite{url:LeadSheet}.
There are some interesting differences, however, which serve to enhance the harmonic structure in the
recorded version.

To transcribe the chords in the recorded version of \IML, we used a computer program called \Chordino.
\Chordino\ is available from the link given in \cite{url:Chordino}. 
We used \Chordino\ as a plug-in to the free audio software \Audacity. To analyze the recording of
\IML, we first converted the stereo file from the \emph{Rubber Soul} CD to
mono.\footnote{We found that the lossless *.wav form of the recording did not produce good
results, as \Chordino\ could not reliably reproduce known chords from the lead sheet. However,
when we analyzed *.m4a (highest quality of $320$~kbs) and *.mp3 (standard quality $128$~kbs) formats,
\Chordino\ reproduced the chords found in the lead sheet and gave the same results for additional
chords, and these are the chords we report here.} We used
the default settings in \Chordino, with one modification: Following the advice given
in \cite{url:Chordino} for analyzing popular music, we set the spectral roll-off parameter to $1.0$. \Chordino\ identifies
chords by detecting the harmonics for the notes composing the chords. In
Figures~\ref{fig:VoiceLeadingsAndHarmonicsInVerseClip} and \ref{fig:VoiceLeadingAndHarmonicsInBridgeClip}
you can gain some idea how idenfifying the sets of harmonics in the sound can be used to identify
the individual notes in a chord, and thus identify the chords themselves.

The chordal analysis given by \Chordino\ differed from the lead sheet in just a few places. The
differences occurred in the $\SpecialText{Verse}$ and $\SpecialText{Bridge}$ parts. Here we will
discuss the differences for the first instance of $\SpecialText{Verse}$ and the second
instance of $\SpecialText{Bridge}$. That will capture the essence of the differences between
the lead sheet version and the recording.
\begin{figure}[!htb]\centering
\setlength{\unitlength}{1in}
\begin{picture}(6.3,2.3)
\put(0,0){\resizebox{6.3in}{!}{\includegraphics{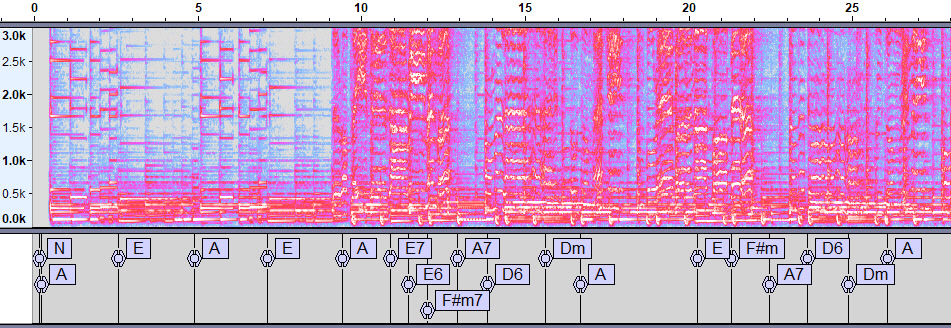}}}
\end{picture}
\FigureCaption{\Chordino\ analysis of chords in first instance of $\SpecialTextFN{Verse}$. The notation $\mathrm{N}$ indicates that no identifiable chord is played at the start, so the first chord is $\mathrm{A}$.\label{fig:ChordinoAnalysisVerse}}
\end{figure}
\subsection{Harmonic Structure of $\SpecialText{Verse}$\label{subsec:HarmonicStructureVerseRecording}}%
In Figure~\ref{fig:ChordinoAnalysisVerse}, we show the \Chordino\ analysis of the chords in
the first instance of $\SpecialText{Verse}$.
This analysis gives us the following progressions:
\begin{equation}\label{eq:ChordinoProgressionsInVerse}
\begin{split}
&\TonnMap{A}{E}\TonnMapB{A}\TonnMapB{E}\TonnMapB{A}\TonnMapB{E}
\TonnMapTTB{f^\sharp}\TonnMapB{A}\TonnMapB{D}\\[3pt]
&\TonnMapB{d}\TonnMapB{A}\TonnMapB{E}\TonnMapTTB{f^\sharp}\TonnMapB{A}\TonnMapB{D}\TonnMapB{d}\TonnMapB{A}
\end{split}
\end{equation}
where we have replaced some chords by their underlying triads---viewing $\mathrm{E^7}, \mathrm{E^6}$ as embellishments of
the chord $\mathrm{E}$, $\mathrm{f^{\sharp 7}}$ as an embellishment of the chord $\mathrm{f^\sharp}$, $\mathrm{D^6}$ as an embellishment of the chord $\mathrm{D}$, and $\mathrm{A^{\! 7}}$ as an embellishment of the chord $\mathrm{A}$.

\begin{figure}[!htb]
\setlength{\unitlength}{1in}
\begin{picture}(4,3.45)\thicklines
\put(0.6,3.3){\SpecialText{Verse}}
\put(0,0){\resizebox{1.5in}{!}{\includegraphics{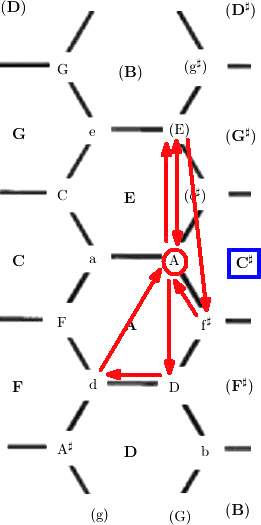}}}
\put(2.9, 3.3){\SpecialText{2nd Bridge}}
\put(2.5,0){\resizebox{1.5in}{!}{\includegraphics{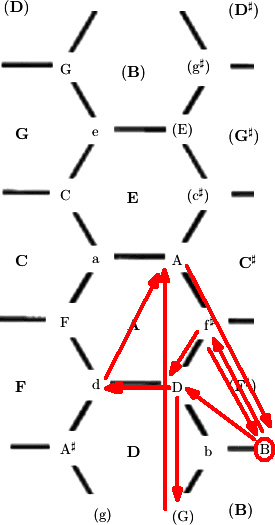}}}
\end{picture}
\FigureCaption{Movement of chords on the \Tonnetz\ for the first instance of $\SpecialTextFN{Verse}$ and the second instance of $\SpecialTextFN{Bridge}$. The blue rectangle in the $\SpecialTextFN{Verse}$ diagram marks the note $\mathrm{C^\sharp}$ that immediately precedes the chord $\mathrm{f^\sharp}$ in the progression $\TonnMapTT{E}{f^\sharp}$.\label{fig:TonnetzDiagramsForInMyLifeChordinoAnalysis}}
\end{figure}

Unlike in the lead sheet $\SpecialText{Verse}$, there are a couple of double \Tonnetz\ transformations in addition
to the many single \Tonnetz\ transformations. The motion of these chord progressions on the \Tonnetz\ is shown
on the left of Figure~\ref{fig:TonnetzDiagramsForInMyLifeChordinoAnalysis}.
The overall pattern is quite similar to the
one shown for $\SpecialText{Verse}$ in Figure~\ref{fig:TonnetzDiagramsForInMyLife}. The one difference is the
two double \Tonnetz\ transformations of type $\TonnMapTT{E}{f^\sharp}$.  These double \Tonnetz\ transformations
enhance the pattern from the lead sheet. To be more precise, in the key of $\mathrm{A}$, they are both
instances of the deceptive progression: $\TonnMap{\mathbf{V}}{\mathbf{vi}}$, which adds some variety to the
harmony. Moreover, the occurrence of the notes $\mathrm{C^\sharp}$ in Lennon's lyrics immediately preceding the chords $\mathrm{f^\sharp}$ makes the transformations $\TonnMapTT{E}{f^\sharp}$ more sonorous.

\begin{figure}[!b]\centering
\setlength{\unitlength}{1.2in} 
\resizebox{6.79in}{!}{%
\begin{picture}(6.5,1.5)
\put(0,0){
\begin{picture}(6.5,1.5)
\put(0,0){
\begin{picture}(1.5,1.5)
\put(0.56,1.42){$\mathrm{A}$}\put(0.53,-0.02){$\mathrm{E}$}
\put(0,0.15){\resizebox{1.44in}{!}{\includegraphics{InMyLifeWholeNoteHarmonicRhythmClock.jpg}}}
\end{picture}}
\put(1.3,1.0){\vector(1,0){0.3}}
\put(1.6,0){
\begin{picture}(1.5,1.5)
\put(-0.15,0.71){$\mathrm{A}$}\put(0.56,1.42){$\mathrm{A}$}\put(1.25,0.7){$\mathrm{E}$}\put(0.53,-0.02){$\mathrm{f^\sharp}$}
\put(0,0.15){\resizebox{1.44in}{!}{\includegraphics{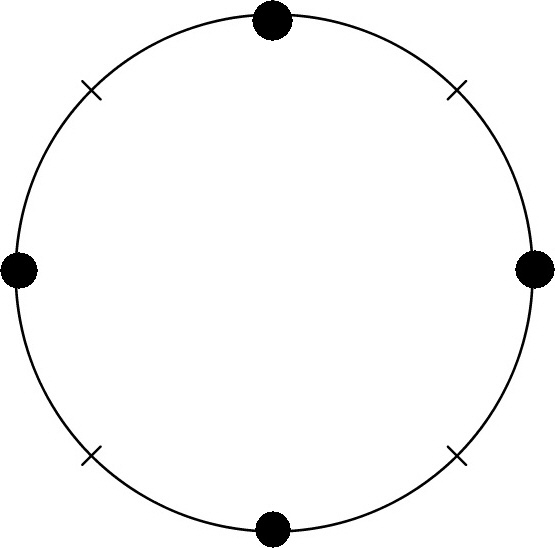}}}
\end{picture}}
\put(2.9,1.0){\vector(1,0){0.3}}
\put(3.2,0){
\begin{picture}(1.5,1.5)
\put(0.53,1.42){$\mathrm{D}$}\put(1.25,0.7){$\mathrm{d}$}\put(0.53,-0.02){$\mathrm{A}$}
\put(0,0.15){\resizebox{1.44in}{!}{\includegraphics{InMyLifeHalfNoteHarmonicRhythmClockVersion2.jpg}}}
\end{picture}}
\put(4.5,1.0){\vector(1,0){0.3}}\put(4.7,1.0){\vector(-1,0){0.2}}
\put(4.8,0){
\begin{picture}(1.5,1.5)
\put(-0.15,0.71){$\mathrm{A}$}
\put(1.25,0.7){$\mathrm{E}$}\put(0.53,-0.02){$\mathrm{f^\sharp}$}
\put(0,0.15){\resizebox{1.44in}{!}{\includegraphics{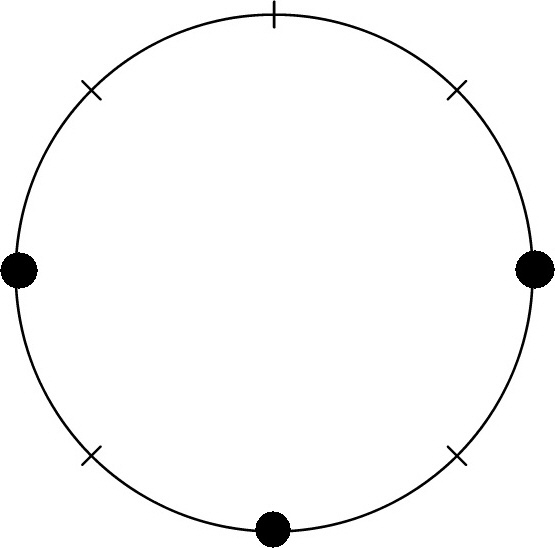}}}
\end{picture}}
\end{picture}}
\end{picture}}
\FigureCaption{Harmonic rhythm of first version of $\SpecialTextFN{Verse}$ in the recording of \IML. The meaning of the hours is the same as for the clocks in Figure~\ref{fig:HarmonicRhythmOfVerse}. The first clock is cycled through twice. The second clock is cycled through once. The third and fourth clocks are cycled through in alternation (third clock, fourth clock, third clock). \label{fig:HarmonicRhythmOfVerseRecordedVersion}}
\end{figure}

The sequence in \eqref{eq:ChordinoProgressionsInVerse} breaks up into four distinct substructures:
\begin{equation}\label{eq:RecordedVerseSubstructures}
\TonnMap{A}{E}\qquad\quad \TonnMap{A}{E}\TonnMapTTB{f^\sharp}\RiemannMapB{R}{A}\qquad\quad
\TonnMap{A}{D}\RiemannMapB{P}{d}\qquad\quad\TonnMap{A}{E}\TonnMapTTB{f^\sharp}
\end{equation}
which differ from the three substructures in the lead sheet $\SpecialText{Verse}$ given in \eqref{eq:SubstructuresForVerseProgressions}. Likewise, there is a different harmonic rhythm corresponding to
these substructures. This harmonic rhythm is shown in Figure~\ref{fig:HarmonicRhythmOfVerseRecordedVersion}.

Thus, in the recorded version of \IML\ there is a slightly more sophisticated sequence of chord progressions (including
enhanced voice leading), as well as slightly more sophisticated harmonic rhythm. This increased sophistication also
occurs in the $\SpecialText{Bridge}$ part.
\subsection{Harmonic Structure of $\SpecialText{Bridge}$\label{subsec:HarmonicStructureOfBridgeInRecording}}
\begin{figure}[!htb]\centering
\setlength{\unitlength}{1in}
\begin{picture}(5.5,2.4)
\put(0,0){\resizebox{5.5in}{!}{\includegraphics{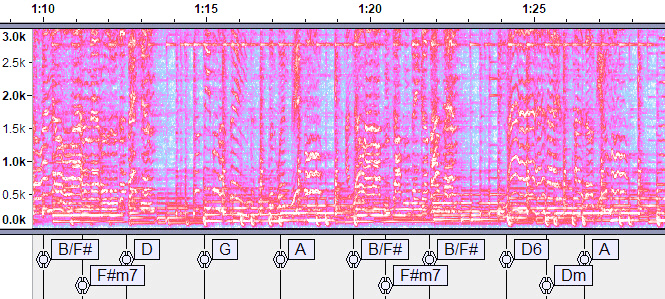}}}
\end{picture}
\FigureCaption{\Chordino\ analysis of chords in second version of $\SpecialTextFN{Bridge}$. \label{fig:ChordinoAnalysisBridgeInterlude}}
\end{figure}

The second instance of $\SpecialText{Bridge}$ most clearly illustrates the enhanced harmonic structure
of the recorded version of \IML. In Figure~\ref{fig:ChordinoAnalysisBridgeInterlude}, we show the
\Chordino\ analysis of the chords in the second instance of $\SpecialText{Bridge}$.
The chord progressions identified by \Chordino\ are the following:
\begin{equation}\label{eq:ChordinoProgressionsInBridge}
\TonnMap{B}{f^\sharp}\TonnMapB{D}\TonnMapB{G}\TonnMapTTB{A}\TonnMapTTB{B}\TonnMapB{f^\sharp}
\TonnMapB{B}\TonnMapB{D}\TonnMapB{d}\TonnMapB{A}
\end{equation}
where we have written the chord $\mathrm{B/F^\sharp}$ as simply $\mathrm{B}$ and written a couple of other chords by
their underlying triads---viewing $\mathrm{f^{\sharp 7}}$ as an embellishment of the chord $\mathrm{f^\sharp}$, and $\mathrm{D^6}$ as an embellishment of the chord $\mathrm{D}$.

As with the lead sheet version, most of these transformations are single \Tonnetz\ transformations. The pattern of the
progressions on the \Tonnetz\ is shown on the right of Figure~\ref{fig:TonnetzDiagramsForInMyLifeChordinoAnalysis}.
As with $\SpecialText{Verse}$, this second version of $\SpecialText{Bridge}$ is quite similar to its lead sheet
version, but with some interesting enhancements. For one thing, the tonicization of the subdominant $\mathrm{D}$ chord
is more clearly emphasized as the central chord in the diagram in Figure~\ref{fig:TonnetzDiagramsForInMyLifeChordinoAnalysis}.
For another, the chromatic chord $\mathrm{B}$ receives much more emphasis. This enhanced emphasis includes the
use of the chromatic mediant transformation $\TonnMap{B}{D}$.

The harmonic rhythm in the recorded version is also more interesting than in the lead sheet. In the recorded version,
if we look at the length of separation between successive chords shown in Figure~\ref{fig:ChordinoAnalysisBridgeInterlude},
we can see that this rhythm consists of interspersing half-note rhythms (such as $\mathrm{B/F^\sharp}$ to $f^{\sharp 7}$) and
whole note rhythms (such as $\mathrm{D}$ to $\mathrm{G}$) with one slightly different splitting of one measure
($\mathrm{B/F^\sharp}$ to $f^{\sharp 7}$ to $\mathrm{B/F^\sharp}$). This is a much more nuanced harmonic rhythm than
in the lead sheet.
\subsection{$\SpecialText{Interlude}$ and $\SpecialText{Coda}$\label{subsec:InterludeAndCodaRecordedVersion}}%
The $\SpecialText{Interlude}$ and $\SpecialText{Coda}$ of the recorded version are identical
in their harmonic structure to the lead sheet version. This is understandable as the $\SpecialText{Interlude}$
was composed for solo piano (played by George Martin and then sped up to sound like a harpsichord). Consequently,
the lead sheet is a faithful transcription of the music (undoubtedly played by Martin from his own sheet music).
The $\SpecialText{Coda}$ is so brief that it probably allowed no room for any elaborations by The Beatles during the
recording process.
\subsection{Summary\label{subsec:SummaryRecording}}%
The recorded version of \IML\ is an enhancement of the lead sheet version. We have
seen how the recorded version builds upon the underlying logic of the lead sheet version.
It does this by an increased emphasis on chromaticism, slightly more complex chord progressions,
and slightly more complex harmonic rhythm.
\section{Conclusion\label{sec:Conclusion}}
In this paper we have used geometric methods to analyze the harmonic progressions and the harmonic
rhythms in the song, \emph{In My Life}. The \Tonnetz\ provided a powerful geometric tool for understanding
the logic of the song's harmony. Spectrograms provided a useful tool for visualizing the interplay of harmonics
in voice leading. Clock diagrams were also quite useful for
analyzing the details of the changing harmonic rhythms in the song. We have endeavored to show as well
that these two geometric approaches were closely linked in the $\SpecialText{Verse}$ parts of \emph{In My Life}.
Our analysis has revealed many of the details of how The Beatles elaborated on basic harmonic patterns to
create progressions more elaborate that typical ones used in rock music at that time.

\end{document}